\documentclass[12pt]{article}
\usepackage{mathrsfs}
\newtheorem{theorem}{Theorem}

\newtheorem{lemma}[theorem]{Lemma}

\usepackage{amssymb}
\usepackage{amsmath}
\usepackage{amsfonts}
\usepackage{graphicx,psfrag}
\usepackage{tikz}
\begin{document}
\def \la{\lambda}

\title{On oriented graphs with minimal skew energy}

\author{Shicai Gong$^{a}$\thanks{Corresponding author. E-mail addresses:
scgong@zafu.edu.cn(S. Gong); lxl@nankai.edu.cn(X. Li);
ghxu@zafu.edu.cn(G. Xu).} \thanks{ Supported by  Zhejiang Provincial
Natural Science Foundation of China(No. Y12A010049).}, Xueliang
Li$^{b}$\thanks{ Supported by National Natural Science Foundation of
China(No. 10831001).} $~$ and Guanghui Xu$^{a}$\thanks{ Supported by
National Natural Science Foundation of China(No. 11171373).}\\
\\{\small \it a. Zhejiang A $\&$ F University, Hangzhou, 311300, P.
R. China}
\\{\small  \it b.Center for Combinatorics and LPMC-TJKLC, Nankai
University,}\\ {\small  \it Tianjin 300071, P. R. China}
   }
\date{}
\maketitle

\begin{abstract}
Let $S(G^\sigma)$ be the skew-adjacency matrix of an oriented graph
$G^\sigma$. The skew energy of $G^\sigma$ is defined as the sum of
all singular values of its skew-adjacency matrix $S(G^\sigma)$. In
this paper, we first deduce an integral formula for the skew energy
of an oriented graph. Then we determine all oriented graphs with
minimal skew energy among all connected oriented graphs on $n$
vertices with $m \ (n\le m < 2(n-2))$ arcs, which is an analogy to
the conjecture for the energy of undirected graphs proposed by
Caporossi {\it et al.} [G. Caporossi, D. Cvetkovi$\acute{c}$, I.
Gutman, P. Hansen, Variable neighborhood search for extremal graphs.
2. Finding graphs with external energy, J. Chem. Inf. Comput. Sci.
39 (1999) 984-996.].

\vskip 0.3cm


\end{abstract}

\section{Introduction}

Let $G^\sigma$ be a digraph that arises from a simple undirected
graph $G$ with an orientation $\sigma$, which assigns to each edge
of $G$ a direction so that $G^\sigma$ becomes an \emph{oriented
graph}, or a \emph{directed graph}. Then $G$ is called the
\emph{underlying graph} of $G^\sigma$. Let $G^\sigma$ be an
undirected graph with vertex set
$V(G^\sigma)=\{v_1,v_2,\cdots,v_n\}.$ Denote by $(u, v)$ an arc, of
$G^\sigma$, with tail $u$ and head $v$. The {\it skew-adjacency
matrix} related to $G^\sigma$ is the $n \times n$ matrix
$S(G^\sigma) = [s_{ij} ],$ where the $(i,j)$ entry satisfies:
$$s_{ij}=\left \{ \begin{array}{ll}
1, & {\rm if \mbox{ } (v_i, v_j)\in G^\sigma \mbox{ }};\\
-1, & {\rm if \mbox{ } (v_j, v_i)\in G^\sigma \mbox{ } };\\
0, & {\rm otherwise. \mbox{ }}
\end{array}\right.$$
The {\it skew energy} of an oriented graph $G^\sigma$, introduced by
Adiga, Balakrishnan and So in \cite{ad} and denoted by
$\mathcal{E}_S(G^\sigma)$, is defined as the sum of all singular
values of $S(G^\sigma)$. Because the skew-adjacency matrix
$S(G^\sigma)$ is skew-symmetric,  the eigenvalues $\{\la_1, \la_2,
\cdots, \la_n\}$ of $S(G^\sigma)$ are all purely imaginary numbers.
Consequently, the skew energy $\mathcal{E}_S(G^\sigma)$ is  the sum
of the absolute values of its eigenvalues, {\it
i.e.,}$$\mathcal{E}_S(G^\sigma)=\sum_{i=1}^n|\la_i|,$$which has the
same expression as that of the energy of an undirected graph with
respect to its adjacent matrix; see e.g. \cite{iv}.

The work on the energy of a graph can be traced back to 1970's
\cite{g1} when Gutman investigated the energy with respect to the
adjacency matrix of an undirected graph, which has a still older
chemical origin; see e.g. \cite{cou}. Then much attention has been
devoted to
 the energy of the adjacency matrix of a graph; see e.g.
 \cite{aa,ago,bs,ds,gkm,gkmz,iv1,yp3,ls,lz1,zz}, and the references
cited therein.  For undirected graphs,  Caporossi,
Cvetkovi$\acute{c}$, Gutman and Hansen \cite{cc} proposed a
conjecture for the minimum energy    as follows.

{\bf Conjecture 1.} Let $G$ be the graph with minimum energy among
all connected graphs  with $n\ge 6$ vertices and $m \ (n\le m \le
2(n-2))$ edges. Then $G$ is $O_{n,m}$ if $m\le n+\lfloor
\frac{n-7}{2}\rfloor$; and $B_{n,m}$ otherwise, where $O_{n,m}$ and
$B_{n,m}$ are respectively the underlying graphs of the oriented
graphs $O^+_{n,m}$ and $B^+_{n,m}$ given in Fig. 1.1.

This conjecture was proved to be true for $m=n-1, 2(n-2)$ by
Caporossi et al. (\cite{cc}, Theorem 1), and $m=n$ by Hou
\cite{yp1}. In \cite{lz1}, Li, Zhang and Wang confirmed this
conjecture for bipartite graphs. Conjecture 1 has not yet been
solved completely.

Recently, in analogy to the energy of the adjacency matrix, a few
other versions of graph energy were introduced in the mathematical
literature, such as Laplacian energy \cite{gz}, signless Laplacian
energy \cite{gr} and skew energy \cite{ad}.

In \cite{ad}, Adiga {\it et.al.} obtained the skew energies of
directed cycles under different orientations and showed that the
skew energy of a directed tree is independent of its orientation,
which is equal to the energy of its underlying tree. Naturally, the
following question is interesting:

\noindent {\bf Question:} Denote by $M$ a class of oriented graphs.
Which oriented graphs have the extremely skew energy among all
oriented graphs of $M$ ?

Hou et al. \cite{hou1} determined the oriented unicyclic graphs with
the maximal and minimal skew energies. Zhu \cite{Zhu} determined the
oriented unicyclic graphs with the first $\lfloor \frac {n-9}
2\rfloor$ largest skew energies. Shen el al. \cite{Shen} determined
the bicyclic graphs with the maximal and minimal energies. Gong and
Xu \cite{GX} determined the 3-regular graphs with the optimum skew
energy, and Tian \cite{Tian} determined the hypercubes with the
optimum skew energy. In the following we will study the minimal skew
energy graphs of order $n$ and size $m$.

At first, we need some notations. Denote by $K_n$, $S_n$ and $C_n$
the complete undirected graph, the undirected star and the
undirected cycle on $n$ vertices, respectively. Let $ O^+_{n,m}$ be
the oriented graph on $n$ vertices which is obtained from the
oriented star $S^\sigma_n$
 by adding $m-n+1$ arcs such that all those arcs have a common vertex; see Fig. 1.1, where $v_1$ is the tail of each arc
incident to it and $v_2$ is the head of each arc incident to it, and
$ B^+_{n,m}$, the oriented graph  obtained from $ O^+_{n,m+1}$ by
deleting the arc $(v_1,v_2)$. Denote by $ O_{n,m}$ and $ B_{n,m}$
the underlying graphs of $ O^+_{n,m}$ and $ B^+_{n,m}$,
respectively. Notice that  both $ O^+_{n,m}$ and $ B^+_{n,m}$
contain $n$ vertices and $m$ arcs.

 \setlength{\unitlength}{0.8mm}
\begin{center}
\begin{picture}(140,31)
\put(10,20){\circle*{2}}  
\put(10,14){\circle*{2}} \put(10,0){\circle*{2}}
\put(10,6){$\vdots$}

\put(20,10){\circle*{2}}  \put(18,5){$v_1$}

\put(19,11){\vector(-1,1){8}} \put(19,10){\vector(-2,1){8}}
\put(19,9){\vector(-1,-1){8}}

\put(40,10){\circle*{2}} \put(38,5){$v_2$}

\put(17,-15){$O^+_{n,m}$}

\put(21,10){\vector(1,0){17}} \put(21,10){\line(1,0){18}}

\put(30,22){\circle*{2}} \put(30,14.5){\circle*{2}}
\put(29.5,16){$\vdots$}

\put(21,10.5){\vector(2,1){8}} \put(21,10.5){\vector(3,4){8}}

\put(39,10.5){\line(-2,1){8}} \put(39,10.5){\line(-3,4){8}}
\put(30,15){\vector(2,-1){7}} \put(30,22.5){\vector(3,-4){7}}

\put(30,1){\circle*{2}}

\put(21,9.5){\vector(1,-1){7}}\put(21,9.5){\line(1,-1){9}}

\put(30,.5){\vector(1,1){7}}\put(39,9.5){\line(-1,-1){9}}

\put(90,20){\circle*{2}}

\put(90,14){\circle*{2}} \put(90,0){\circle*{2}}
\put(90,6){$\vdots$}

\put(100,10){\circle*{2}}  \put(98,5){$v_1$}

\put(99,11){\vector(-1,1){8}} \put(99,10){\vector(-2,1){8}}
\put(99,9){\vector(-1,-1){8}}

\put(120,10){\circle*{2}} \put(118,5){$v_2$}

\put(97,-15){$B^+_{n,m}$}

\put(110,22){\circle*{2}}

\put(110,14.5){\circle*{2}}

\put(109.5,16){$\vdots$}

\put(101,10.5){\vector(2,1){8}} \put(101,10.5){\vector(3,4){8}}

\put(119,10.5){\line(-2,1){8}} \put(119,10.5){\line(-3,4){8}}
\put(110,15){\vector(2,-1){7}} \put(110,22.5){\vector(3,-4){7}}

\put(110,1){\circle*{2}}

\put(101,9.5){\vector(1,-1){7}}\put(101,9.5){\line(1,-1){9}}

\put(110,.5){\vector(1,1){7}}\put(119,9.5){\line(-1,-1){9}}

\put(16,-26) {Fig. 1.1. Two oriented graphs $O^+_{n,m}$ and
$B^+_{n,m}$. }

\end{picture}
\end{center}
\vspace{2.5cm}

In this paper, we  first deduce an integral formula for the skew
energy of an oriented graph. Then we study  the question above and
determine all oriented graphs with minimal skew energy among all
connected oriented digraphs on $n$ vertices with $m \ (n\le m <
2(n-2))$ arcs. Interestingly, our result is an analogy to Conjecture
1.
\begin{theorem} \label{01} Let
$G^\sigma$ be an oriented graph with minimal skew energy among all
oriented  graphs with $n$ vertices and $m \ (n\le m < 2(n-2))$ arcs.
Then, up to isomorphism, $G^\sigma$ is\\ $(1)$ $O^+_{n,m}$ if
$m<\frac{3n-5}{2}$; \\ $(2)$ either $B^+_{n,m}$ or $O^+_{n,m}$ if
$m=\frac{3n-5}{2}$; and \\ $(3)$ $B^+_{n,m}$ otherwise.
\end{theorem}

\vskip 0.5cm

\section{Integral formula for the skew energy}

In this section, based on the formula established by Adiga {\it et
al.} \cite{ad}, we deduce an integral formula for the skew energy of
an oriented graph, which is an analogy to the Coulson integral
formula for the energy of an undirected graph. Firstly, we introduce
some notations and preliminary results.

An even cycle $C$ in an oriented graph $G^\sigma$ is called
\emph{oddly oriented} if for either choice of direction of traversal
around $C$, the number of edges of $C$ directed in the direction of
the traversal is odd. Since $C$ is even, this is clearly independent
of the initial choice of direction of traversal. Otherwise, such an
even cycle $C$ is called as \emph{evenly oriented}. (Here we do not
involve the parity of the cycle with length odd. The reason is that
it depends on the initial choice of direction of traversal.)

A ``\emph{basic oriented graph}" is an oriented graph whose
components are even cycles and/or complete oriented graphs with
exactly two vertices.

Denote by $\phi(G^\sigma;x)$ the \emph{ skew characteristic
polynomial} of an oriented graph $G^\sigma$, which is defined as
\begin{eqnarray*}
   \phi(G^\sigma;x)=det(xI_n-S(G^\sigma))=\sum_{i=0}^na_{i}(G^\sigma)x^{n-i},
   \end{eqnarray*}
where $I_n$ denotes the identity matrix of order $n$. The following
result is a cornerstone of our discussion below, which determines
all coefficients of the skew characteristic polynomial of an
oriented graph in terms of its basic oriented subgraphs; see
\cite[Theorem 2.4]{hou} for an independent version.

\begin{lemma}  \em{\cite[Corollary 2.3]{gx}} \label{0001}
Let $G^\sigma$ be an oriented graph on $n$ vertices, and let the
skew characteristic polynomial of $G^\sigma$ be
$$\phi(G^\sigma,\la)=\sum_{i=0}^n(-1)^ia_i
\la^{n-i}=\la^n-a_1\la^{n-1}+a_2\la^{n-2}+\cdots+(-1)^{n-1}a_{n-1}\la+(-1)^na_n.$$
Then $a_i=0$ if $i$ is odd; and
$$a_i=\sum_{\mathscr{H}}(-1)^{c^+}2^c  {\rm \mbox{ } \mbox{ } if \mbox{ }}
i {\rm \mbox{ } is \mbox{ } even},$$ where the summation is over
all basic oriented subgraphs $\mathscr{H}$ of $G^\sigma$ having $i$
vertices  and $c^{+}$ and $c$ are respectively the number of evenly
oriented even cycles and   even cycles contained in $\mathscr{H}$.
\end{lemma}

Let $G=(V(G),E(G))$ be a graph, directed or not, on $n$ vertices.
Then denote by $\Delta(G)$  the maximum degree of $G$ and set
$\Delta(G^\sigma)=\Delta(G)$. An \emph{$r$-matching} in a graph $G$
is a subset of $r$ edges such that every vertex of $V(G)$ is
incident with at most one edge in it. Denote by $M(G,r)$ the number
of all $r$-matchings in $G$ and set $M(G,0)=1$.

Denote by $q(G)$ the number of quadrangles in a undirected graph
$G$. Then as a consequence of Lemma \ref{0001}, we have

\begin{theorem} \label{013} Let $G^\sigma$ be an
oriented graph containing $n$ vertices and $m$ arcs. Suppose
$$\phi(G^\sigma,\la)=\sum_{i=0}^{n}(-1)^ia_{i}(G^\sigma) \la^{n-i}.$$
Then $a_0(G^\sigma)=1$, $a_2(G^\sigma)=m$ and $a_4(G^\sigma)\ge
M(G,2)-2q(G)$ with equality if and only if all oriented quadrangles
of $G^\sigma$ are evenly oriented.
\end{theorem}

\noindent {\bf Proof.} The result follows from Lemma \ref{0001} and
the fact that each arc corresponds a basic oriented graph having $2$
vertices, and each basic oriented graph having $4$ vertices is
either a $2$-matching or a  quadrangle.\hfill $\blacksquare$

Furthermore, as well-known, the eigenvalues of an arbitrary real
skew symmetric matrix  are all purely imaginary numbers and
 occur in conjugate pairs. Henceforth, Lemma \ref{0001} can be
 strengthened as follows,  which will provide  much  convenience for our
discussion below.

\begin{lemma}  \label{2} Let $G^\sigma$ be an
oriented graph of order $n$. Then each coefficient of the skew
characteristic polynomial $$\phi(G^\sigma,\la)=\sum_{i=0}^{\lfloor
\frac{n}{2} \rfloor}a_{2i}(G^\sigma) \la^{n-2i}$$ satisfies
$a_{2i}(G^\sigma)\ge 0$ for each $i(0\le i\le \lfloor \frac{n}{2}
\rfloor)$.
\end{lemma}
{\bf Proof.} Let $\la_1,\la_2,\cdots,\la_n $ be all eigenvalues of
the skew adjacency matrix $S(G^\sigma)$ of $G^\sigma$. Because
$\la_1,\la_2,\cdots,\la_n $ are all purely imaginary numbers and
must occur in conjugate pairs, we suppose, without loss of
generality, that there exists an integer number $m(\le \lfloor
\frac{n}{2} \rfloor)$ such that
$$\la_t=-\la_{n-t+1}=p_t i,\ \ for \ \  t=1,2,\cdots,m,$$ and all other
eigenvalues are zero, where each $p_t$ is a positive real number and
$i$ satisfies $i^2=-1.$ Then we have \begin{equation*}
\begin{array}{lll}
 \phi(G^\sigma,\la)
 &=&\prod_{t=1}^n(\la-\la_t)\\&=&\la^{n-2m}\prod_{t=1}^m(\la^2+p^2_t),
\end{array}
\end{equation*}
which implies that the result follows. \hfill $\blacksquare$
\vspace{3mm}

For  an oriented graph $G^\sigma$ on $n$ vertices, an integral
formula for the skew energy in terms of the skew characteristic
polynomial $\phi(G^\sigma, \la)$ and its derivative is given by
\cite{ad}
$$ \mathscr{E}_s(G^\sigma) =\frac{1}{\pi}\int_{-\infty}^{+\infty}\left[n+\la \frac{\phi'(G^\sigma, -\la)}{\phi(G^\sigma,
-\la)}\right]d\la.\eqno{(2.1)}$$ However, using the above integral,
it is by no means easy to calculate the skew energy of an oriented
graph. Hence, it is rather important to establish some other more
simpler formula.

Applying to (2.1) the fact that the coefficient $a_i=0$ for each odd
$i$ from Lemma \ref{0001} and replacing $\la$ by $-\la$, we have
$$ \mathscr{E}_s(G^\sigma)
=\frac{1}{\pi}\int_{-\infty}^{+\infty}\left[n-\la
\frac{\phi'(G^\sigma, \la)}{\phi(G^\sigma, \la)}\right]d\la.$$
Meanwhile, note that
$$\frac{\phi'(G^\sigma,
\la)}{\phi(G^\sigma, \la)}d\la=d\ln \phi(G^\sigma, \la).$$Then we
have
\begin{equation*}
\begin{array}{lll}
 \mathscr{E}_s(G^\sigma)
&=&\frac{1}{\pi}\int_{-\infty}^{+\infty}\left[n-\la
\frac{\phi'(G^\sigma, \la)}{\phi(G^\sigma,
\la)}\right]d\la \\
  &=&\frac{1}{\pi}\int_{-\infty}^{+\infty}\left[n-\la (\frac{d}{d\la}) \ln\phi(G^\sigma,
\la)\right]d\la.
\end{array}\eqno{(2.2)}
\end{equation*}
Therefore, we have
\begin{theorem} \label{04} Let $G^\sigma$ be an
oriented graph with order $n$. Then $$
\mathscr{E}_s(G^\sigma)=\frac{1}{\pi}\int_{-\infty}^{+\infty}\la^{-2}\ln
\psi(G^\sigma,\la)d\la,\eqno{(2.3)}$$where
$$\psi(G^\sigma,\la)=\sum_{i=0}^{\lfloor \frac{n}{2} \rfloor}a_{2i}(G^\sigma) \la^{2i}$$ and $a_{2i}(G^\sigma)$ denotes
the coefficient of $\la^{n-2i}$ in the skew characteristic
polynomial $\phi(G^\sigma,\la)$.
\end{theorem}
{\bf Proof.} Let both $G^{\sigma_1}_1$ and $G^{\sigma_2}_2$ be
oriented graphs with order $n$. ($G_1$ perhaps equals $G_2$.) Then
applying (2.2) we have
\begin{eqnarray*}
 \mathscr{E}_s(G^{\sigma_1}_1)-\mathscr{E}_s(G^{\sigma_2}_2) =-\frac{1}{\pi}\int_{-\infty}^{+\infty}\la
(\frac{d}{d\la}) \ln\left[\frac{\phi(G^{\sigma_1}_1,
\la)}{\phi(G^{\sigma_2}_2, \la)}\right]d\la.
\end{eqnarray*}
Using partial integration, we have
$$\mathscr{E}_s(G^{\sigma_1}_1)-\mathscr{E}_s(G^{\sigma_2}_2)
=-\frac{\la}{\pi}\ln[ \frac{\phi(G^{\sigma_1}_1,
\la)}{\phi(G^{\sigma_2}_2,
\la)}]|_{-\infty}^{+\infty}+\frac{1}{\pi}\int_{-\infty}^{+\infty}\ln\left[
\frac{\phi(G^{\sigma_1}_1, \la)}{\phi(G^{\sigma_2}_2,
\la)}\right]d\la.
$$ Notice that
$$\frac{\la}{\pi}\ln\left[ \frac{\phi(G^{\sigma_1}_1, \la)}{\phi(G^{\sigma_2}_2,
\la)}\right]\mid_{-\infty}^{+\infty}=0.$$ Hence
$$ \mathscr{E}_s(G^{\sigma_1}_1)-\mathscr{E}_s(G^{\sigma_2}_2)=\frac{1}{\pi}\int_{-\infty}^{+\infty}\ln\left[
\frac{\phi(G^{\sigma_1}_1, \la)}{\phi(G^{\sigma_2}_2,
\la)}\right]d\la.$$ Suppose now that $G^{\sigma_2}_2$ is the null
oriented graph, an oriented graph containing $n$ isolated vertices.
Then $\phi(G^{\sigma_2}_2, \la)=\la^n$ and thus
$\mathscr{E}_s(G^{\sigma_2}_2)=0$. After an appropriate change of
variables we can derive
$$\mathscr{E}_s(G^{\sigma_1}_1)=\frac{1}{\pi}\int_{-\infty}^{+\infty}\la^{-2}\ln
\psi(G^{\sigma_1}_1,\la)d\la.$$ Then the result follows.\hfill
$\blacksquare$

\section{Proof of Theorem \ref{01}}

From  Theorem \ref{04}, for an oriented graph $G^\sigma$ on $n$
vertices, the skew energy $\mathcal{E}_s(G^\sigma)$ is a strictly
monotonically increasing function of the coefficients
$a_{2k}(G^\sigma) (k = 0, 1, \cdots, \lfloor\frac{n}{2}\rfloor)$,
since for each $i$ the coefficient of $\la^{n-i}$ in the
characteristic polynomial $\phi(G^\sigma,\la)$, as well as
$\psi(G^\sigma,\la)$, satisfies $a_i(G^\sigma)\ge 0$ by Lemma
\ref{2}. Thus, similar to comparing two undirected graphs with
respect to their energies, we define the quasi-ordering relation
$``\preceq"$ of two oriented graphs with respect to their skew
energies as follows.\vspace{2mm}

Let $G^{\sigma_1}_1$ and $G^{\sigma_2}_2$ be two oriented graphs of
order $n$. ($G_1$ is not necessarily different from $G_2$.) If
$a_{2i}(G^{\sigma_1}_1)\le a_{2i}(G^{\sigma_2}_2)$ for all $i$ with
$0\le i \le \lfloor\frac{n}{2}\rfloor$, then we write that
$G^{\sigma_1}_1 \preceq G^{\sigma_2}_2$.\vspace{2mm}

Furthermore, if $G^{\sigma_1}_1 \preceq G^{\sigma_2}_2$ and there
exists at least one index $i$ such that $a_{2i}(G^{\sigma_1}_1)<
a_{2i}(G^{\sigma_2}_2)$, then we write that $G^{\sigma_1}_1 \prec
G^{\sigma_2}_2$. If $a_{2i}(G^{\sigma_1}_1)= a_{2i}(G^{\sigma_2}_2)$
for all $i$, we write $G^{\sigma_1}_1 \sim G^{\sigma_2}_2$. Note
that there are non-isomorphic oriented graphs $G^{\sigma_1}_1$ and
$G^{\sigma_2}_2$ such that $G^{\sigma_1}_1 \sim G^{\sigma_2}_2$,
which implies that $``\preceq"$ is not a partial ordering in
general. \vspace{2mm}

According to the integral formula (2.2), we have, for two oriented
graphs $D_1$ and $D_2$ of order $n$, that $$D_1 \preceq
D_2\Longrightarrow \mathcal{E}_s(D_1)\le \mathcal{E}_s(D_2)$$ and
$$D_1 \prec D_2\Longrightarrow \mathcal{E}_s(D_1)<
\mathcal{E}_s(D_2).\eqno{(3.1)}$$

In the following, by discussing the relation ``$\succeq$", we
compare the skew energies for two oriented graphs and then complete
the proof of Theorem \ref{01}.

Firstly, by a directly calculation we have
$$\phi(O^+_{n,m})=\la^n+m\la^{n-2}+(m-n+1)(2n-m-3)\la^{n-4},\eqno{(3.2)}$$and
$$\phi(B^+_{n,m})=\la^n+m\la^{n-2}+(m-n+2)(2n-m-4)\la^{n-4}.\eqno{(3.3)}$$

Denote by $G^\sigma(n,m)$ and $G(n,m)$ the sets of all connected
oriented graphs and undirected graphs with $n$ vertices and $m$
edges, respectively. The following results on undirected graphs are
needed.
\begin{lemma} \label{05} Let $n\ge 5$ and $G\in G(n,m)$ be an arbitrary
connected undirected graph containing $n$ vertices and $m \ (n\le m
<2(n-2))$ edges. Then $q(G)\le \left(
\begin{array}{c}
m-n+2  \\
    2
\end{array}\right)
$, where $q(G)$ denotes the number of quadrangles contained in $G$.
\end{lemma}

{\bf Proof.} We prove this result by induction on $m$.

The result is obvious for $m= n$. So we suppose that $n< m <2(n-2)$
and the result is true for smaller $m$.

Let $e$ be an edge of $G$ and  $q_G(e)$ denote the number of
quadrangles containing the edge $e$. Suppose $e=(u,v)$. Let $U$ be
the set of neighbors of $u$ except $v$, and $V$ the set of neighbors
of $v$ except $u$. Then there are just $q_G(e)$ edges between $U$
and $V$. Let $X$ be the subset of $U$ such that each vertex in $X$
is incident to some of the above $q_G(e)$ edges and $Y$ be the
subset of $V$ defined similarly to $X$. Assume $|X| =x$ and $|Y|=
y$. Let $G_0$ be the subgraph of $G$ induced by $V(G_0)=u \cup v
\cup X \cup Y$. Then there are at least $q_G(e)+x+y+1$ edges and
exactly $x + y +2$ vertices in $G_0$. In order for the remaining
vertices to connect to $G_0$, the number of remaining edges must be
not less than that of the remaining vertices. Thus
$$m-(q_G(e)+x+y+1)\ge n-(x+y+2).$$ That is $$q_G(e)\le m-n+1.$$
By induction hypothesis, $q(G-e)\le \left(\begin{array}{c}
(m-1)-n+2  \\
    2
\end{array}\right)=\left(\begin{array}{c}
m-n+1  \\
    2
\end{array}\right).$
Then we have
$$q(G)=q_G(e)+q(G-e)\le m-n+1+\left(\begin{array}{c}
m-n+1  \\
    2
\end{array}\right)=\left(\begin{array}{c}
m-n+2  \\
    2
\end{array}\right).$$
Hence, the result follows.\hfill $\blacksquare$

By a similar method, we can show that
\begin{lemma} \label{015} Let $n\ge 5$ and $G\in G(n,m)$ be an arbitrary undirected graph
containing $n$ vertices and $m \ (n\le m <2(n-2))$ edges. Suppose
$\Delta(G)=n-1$. Then $$q(G)\le \left(
\begin{array}{c}
m-n+1  \\
    2
\end{array}\right).$$
\end{lemma}
\begin{lemma} {\em \cite[A part of Theorem 2.6]{gx}}\label{016} Let $G^{\sigma}$ be
an oriented graph with an arc $e=(u,v)$. Suppose that $e$ is not
contained in any even cycle. Then
$$\phi(G^{\sigma}, \la)=\phi(G^{\sigma} -e, \la)+ \phi(G^{\sigma}
- u-v,\la).$$
\end{lemma}

As a consequence of Lemma \ref{016}, we have the following result.
\begin{lemma}\label{017} Let $G^{\sigma} $ be an oriented graph on $n$
vertices and $(u,v)$ a pendant arc of $G^{\sigma} $ with pendant
vertex $v$. Suppose $\phi(G^{\sigma},
\la)=\sum_{i=0}^na_i(G^{\sigma})\la^{n-i}.$ Then
$$a_{i}(G^{\sigma})=a_{i}(G^{\sigma}-v)+a_{i-2}(G^{\sigma}-v-u).$$
\end{lemma}

Based on the preliminary results above, we have the following two
results.
\begin{lemma} \label{06} Let $n\ge 5$ and $G^\sigma\in G^\sigma(n,m)$ be
an oriented graph with maximum degree $n-1$. Suppose that $ n \le m
< 2(n-2)$ and $G^\sigma \nsim O^+_{n,m}$. Then $G^\sigma \succ
O^+_{n,m}$.
\end{lemma}
{\bf Proof.} By Theorem \ref{013}, it suffices to prove that
$a_4(G^\sigma)> a_4(O^+_{n,m})$. Suppose that $v$ is the vertex with
degree $n-1$. For convenience, all arcs incident to $v$ are colored
as white and all other arcs are colored as black. Then there are
$n-1$ white arcs and $m-n+1$ black arcs. We estimate the cardinality
of $2$-matchings in $G^\sigma$ as follows. Noticing that all white
arcs are incident to $v$, each pair of white arc can not form a
$2$-matching of $G^\sigma$. Since $d(v)=n-1$ and each black arc
incident to exactly two white arcs, each black arc together with a
white arcs except its neighbors forms a $2$-matching of $G^\sigma$,
that is, there are $(m-n+1)(n-3)$ black-white $2$-matchings.
Moreover, noticing that $G^\sigma\neq O^+_{n,m}$, $G^\sigma-v$ does
not contain the directed star $S_{m-n+2}$ as its subgraph, and thus
there is at least one $2$-matching formed by a pair of disjoint
black arcs, or $G^\sigma$ is an oriented graph of the following
graph $F$. \setlength{\unitlength}{1mm}
\begin{center}
\begin{picture}(50,20)
\put(10,20){\circle*{2}} \put(10,14){\circle*{2}}
\put(10,0){\circle*{2}} \put(10,6){$\vdots$}
\put(20,10){\circle*{2}}  \put(18,5){$v_1$}
\put(19,11){\line(-1,1){9}} \put(19,10){\line(-2,1){9}}
\put(19,9){\line(-1,-1){9}} \put(40,10){\circle*{2}}
\put(38,5){$v_2$}
\put(21,10){\line(1,0){17}}\put(21,10){\line(1,0){18}}
\put(30,18){\circle*{2}} \put(21,10.5){\line(4,3){9}}
\put(39,10.5){\line(-4,3){9}} \put(30,1){\circle*{2}}
\put(21,9.5){\line(1,-1){9}} \put(30,.5){\line(1,1){9}}
\put(30,.5){\line(0,1){18}} \put(8,-9){Fig. 1.2. The graph $F$. }
\end{picture}
\end{center}
\vspace{1cm} If it is the first case, then the number of
$2$-matchings in $G^\sigma$ satisfies
$$M(G^\sigma,2)\ge (m-n+1)(n-3)+1.$$
From Lemma \ref{015}, $q(G^\sigma)\le \left(
\begin{array}{c}
m-n+1  \\
    2
\end{array}\right),$ and then by applying
Theorem \ref{013} again, we have
\begin{eqnarray*} a_4(G^\sigma)&\ge&
M(G^\sigma,2)-2q(G^\sigma)\\&\ge& (m-n+1)(n-3)+1- 2\left(
\begin{array}{c}
m-n+1  \\
    2
\end{array}\right)\\&=&a_4(O^+_{n,m})+1
\end{eqnarray*}
 by Eq.(3.2).
 If it is the second case, clearly $m=n+2$, $q(F)=3$, but the three quadrangles
 can not be all evenly oriented. Then
 \begin{eqnarray*} a_4(F)\ge
M(F,2)-2q(F)\ge (m-n+1)(n-3)-4>a_4(O^+_{n,n+2}).
\end{eqnarray*}
The result thus follows. \hfill $\blacksquare$

\begin{lemma} \label{07} Let $n\ge 5$ and $G^\sigma\in G^\sigma(n,m)$
be an oriented graph with $n\le m < 2(n-2)$. Suppose that
$\Delta(G^\sigma)\le n-2$ and $G^\sigma \nsim B^+_{n,m}$. Then
$G^\sigma \succ B^+_{n,m}$.
\end{lemma}
{\bf Proof.} By Theorem \ref{013} again, it suffices to prove that
$a_4(G^\sigma)> a_4(B^+_{n,m})$. We apply induction on $n$ to prove
it. By a direct calculation, the result follows if $n=5$, since then
$5=m<2(5-2)=6$ and there exists exactly four graphs in
$G^\sigma(5,5)$, namely, the oriented cycle $C_3$ together with two
pendant arcs attached to two different vertices of the $C_3$, the
oddly oriented cycle $C_4$ together with a pendant arc, $B^+_{5,5}$
and the oriented cycle $C_5$. Suppose now that $n\ge 6$ and the
result is true for smaller $n$.

\noindent {\bf Case 1.} There is a pendant arc $(u,v)$ in $G^\sigma$
with pendant vertex $v$.

By Lemma \ref{017} we have
$$a_4(G^\sigma)=a_4(G^\sigma-v)+a_2(G^\sigma-v-u)=a_4(G^\sigma-v)+e(G^\sigma-v-u).$$
Noticing that $\Delta(G^\sigma)\le n-2$, we have $e(G^\sigma-v-u)\ge
m-\Delta(G^\sigma)\ge m-n+2$.

By induction hypothesis, $a_4(G^\sigma-v)\ge a_4(B^+_{n-1,m-1})$
with equality if and only if $G^\sigma-v= B^+_{n-1,m-1}$. Then
\begin{equation*}
\begin{array}{lll}
a_4(G^\sigma)&=&a_4(G^\sigma-v)+a_2(G^\sigma-v-u)\\&\ge&
a_4(B^+_{n-1,m-1})+m-n+2\\&=& a_4(B^+_{n-1,m-1})+e(S_{m-n+1})\\ &=&
a_4(B_{n,m})
\end{array}
\end{equation*}
with equality if and only if $G^\sigma= B^+_{n,m}$. The result thus
follows.

\noindent {\bf Case 2.} There are no pendant vertices in $G^\sigma$.

Let
$$(d)_{G^\sigma}=(d_1,d_2,\cdots,d_i,d_{i+1},\cdots,d_n)$$
be the non-increasing degree sequence of $G^\sigma$. We label the
vertices of $G^\sigma$ corresponding to the degree sequence
$(d)_{G^\sigma}$ as $v_1,v_2,\cdots,v_n$ such that
$d_{G^\sigma}(v_i)=d_i$ for each $i$. Assume $d_1<n-2$. Then there
exists a vertex $v_k$ that is not adjacent to $v_1$, but is adjacent
to one neighbor, say $v_i$, of $v_1$. Thus
$$(d_1+1,d_2,\cdots
d_i-1,d_{i+1},\cdots,d_n)$$ is the degree sequence of the oriented
graph $D'$ obtained from $G^\sigma$ by deleting the arc $(v_k,v_i)$
and adding the arc $(v_k,v_1)$, regardless the orientation of the
arc $(v_k,v_1)$. Rewriting the sequence above such that
$$(d)_{D'}=(d'_1,d'_2,\cdots,d'_{i},d'_{i+1},\cdots,d'_n)$$
is also a non-increasing sequence. Then $d_{1}\ge d_i\ge 2$ and thus
we have
$$\sum_{i=1}^n\left(
\begin{array}{c}
d_i'  \\
    2
\end{array}\right)>
\sum_{i=1}^n\left(
\begin{array}{c}
d_i  \\
    2
\end{array}\right),\eqno{(3.4)}$$since
\begin{equation*}
\begin{array}{lll}\sum_{i=1}^n\left(
\begin{array}{c}
d_i'  \\
    2
\end{array}\right)-
\sum_{i=1}^n\left(
\begin{array}{c}
d_i  \\
    2
\end{array}\right)&=&
\left(
\begin{array}{c}
d_1+1  \\
    2
\end{array}\right)+\left(
\begin{array}{c}
d_i-1  \\
    2
\end{array}\right)-\left(
\begin{array}{c}
d_1  \\
    2
\end{array}\right)-\left(
\begin{array}{c}
d_i  \\
    2
\end{array}\right)\\
&=& d_{1}-d_i+1 \\&>&0.
\end{array}
\end{equation*}
Repeating this procedure, we can eventually obtain a non-increasing
graph sequence
$$(d)_{D''}=(d''_1,d''_2,\cdots,d''_{i},d''_{i+1},\cdots,d''_n)$$such
that $\Delta(D'')=d_1''=n-2$ and
 $$\sum_{v\in D''}\left(
\begin{array}{c}
d''(v)  \\
    2
\end{array}\right)> \sum_{v\in
D' }\left(
\begin{array}{c}
d'(v)  \\
    2
\end{array}\right)>\cdots > \sum_{v\in
G^\sigma}\left(
\begin{array}{c}
d(v)  \\
    2
\end{array}\right).\eqno{(3.5)}$$

Similarly, we can assume that there exists a vertex $v_k$ that is
not adjacent to $v_i$, but is adjacent to one neighbor, say $v_j$,
of $v_i$. Thus
$$(d_1,d_2,\cdots
d_i+1,d_{i+1},\cdots,d_j-1,d_{j+1},\cdots,d_n)$$ is the degree
sequence of the oriented graph $D'''$ obtained from $D''$ by
deleting the arc $(v_k,v_j)$ and adding the arc $(v_k,v_i)$,
regardless the orientation of the arc $(v_k,v_j)$. By a similar
proof, we can get
$$\sum_{v\in D'''}\left(
\begin{array}{c}
d'''(v)  \\
    2
\end{array}\right)> \sum_{v\in
D'' }\left(
\begin{array}{c}
d''(v)  \\
    2
\end{array}\right).$$
Then by applying the above procedure repeatedly, we eventually
obtain the degree sequence $(d)_{B^+_{n,m}}$,
$$(d)_{B^+_{n,m}}=(n-2,m-n+2,2,2,\cdots,2,1,1,\cdots,1),$$
where the number of vertices of degree $2$ is $m-n+2$, and the
number of vertices of degree $1$ is $2n-m-4$. Finally, we get
$$\sum_{v\in B^+_{n,m}}\left(
\begin{array}{c}
d^{B^+}(v)  \\
    2
\end{array}\right)> \sum_{v\in
D''' }\left(
\begin{array}{c}
d'''(v)  \\
    2
\end{array}\right)> \sum_{v\in
D''}\left(
\begin{array}{c}
d''(v)  \\
    2
\end{array}\right)>\cdots > \sum_{v\in
G^\sigma}\left(
\begin{array}{c}
d(v)  \\
    2
\end{array}\right).$$

Then the lemma follows by combining Eq.(3.3) and Lemma \ref{05} with
the fact that $M(G,2)=\left(
\begin{array}{c}
m  \\
    2
\end{array}\right)-\sum_{v\in G^\sigma}\left(
\begin{array}{c}
d(v)  \\
    2
\end{array}\right)$.
\hfill $\blacksquare$

Combining Lemma \ref{06} with Lemma \ref{07}, we get the proof of
Theorem \ref{01} immediately.

\vskip3mm
\noindent {\bf Proof of Theorem \ref{01}.} Combining with
Lemmas \ref{06} and \ref{07}, the oriented graph with minimal skew
energy among all oriented graphs of $G^\sigma(n,m)$ with $n\le m \le
2(n-2)$ is either $O^+_{n,m}$ or $B^+_{n,m}$. Furthermore, from
(3.2) and (3.3), we have
$$a_4(O^+_{n,m})=(m-n+1)(2n-m-3)$$ and
$$a_4(B^+_{n,m})=(m-n+2)(2n-m-4).$$
Then, by a direct calculation we have
$a_4(O^+_{n,m})<a_4(B^+_{n,m})$ if $m<\frac{3n-5}{2}$;
$a_4(B^+_{n,m})=a_4(O^+_{n,m})$ if $m=\frac{3n-5}{2}$; and
$a_4(O^+_{n,m})>a_4(B^+_{n,m})$ otherwise. The proof is thus
complete by (3.1). \hfill $\blacksquare$


\begin{thebibliography}{50}

\bibitem{ad} C. Adiga, R. Balakrishnan and Wasin So, The skew energy of a
graph, \emph{Linear Algebra Appl.} 432 (2010), 1825-1835.

\bibitem{aa} F. Alinaghipour and B. Ahmadi, On the energy of complement
of regular line graph, \emph{MATCH Commun. Math. Comput. Chem.}
60(2008), 427-434.

\bibitem{ago}S. Akbari, E. Ghorbani and M.R. Oboudi, Edge
addition, singular values and energy of graphs and matrices,
\emph{Linear Algebra Appl.} 430(2009): 2192-2199.

\bibitem{bs} S.R. Blackburn and I.E. Shparlinski, On the average energy
of circulant graphs, \emph{Linear Algebra Appl.} 428 (2008),
1956-1963.

\bibitem{cc}G. Caporossi, D. Cvetkovi$\acute{c}$, I. Gutman and P. Hansen, Variable
neighborhood search for extremal graphs. 2. Finding graphs with
external energy, \emph{J. Chem. Inf. Comput. Sci.} 39 (1999),
984-996.

\bibitem{cou} C.A. Coulson, On the calculation of the
energy in unsaturated hydrocarbon molecules, \emph{Proc. Cambridge
Phil. Soc.} 36(1940), 201-203.

\bibitem{ds} J. Day and W. So, Graph energy change due to
edge deletion, \emph{Linear Algebra Appl.}, 428 (2007): 2070-2078.

\bibitem{gx}  S. Gong and G. Xu, The characteristic polynomial and the matchings polynomial
of a weighted oriented graph, \emph{Linear Algebra Appl.} 436(2012),
3597-3607.

\bibitem{GX} S. Gong, G. Xu, 3-Regular digraphs with optimum skew energy,  \emph{Linear Algbra
Appl.} 436(2012), 465-471.

\bibitem{g1} I. Gutman, The energy of a graph, \emph{Ber. Math.-Statist. Sekt.
Forschungsz. Graz.} 103(1978), 1-22.

\bibitem{iv1} I. Gutman, \emph{The energy of a graph: old and new results,}
in: A. Betten, A. Kohnert, R. Laue and A. Wassermann(Eds.),
\emph{Algebraic Combinatorics and Applications,} Springer, Berlin,
2000, pp.196-211.

\bibitem{gkm} I. Gutman, D. Kiani and M. Mirzakhah, On incidence energy of
graphs,\emph{ MATCH Commun. Math. Comput. Chem.} 62 (2009), 573-580.

\bibitem{gkmz} I. Gutman, D. Kiani, M. Mirzakhah and B. Zhou, On incidence energy of a
graph,\emph{ Linear Algebra Appl.} 431 (2009), 1223-1233.

\bibitem{iv} I. Gutman and O.E. Polansky,\emph{ Mathematical Concepts in Organic Chemistry,}
Springer, Berlin, 1986.

\bibitem{gr} I. Gutman, M. Robbiano, E.A. Martins, D.M. Cardoso, L. Medina and O.
Rojo, Energy of line graphs,  \emph{Linear Algebra Appl.} 433(2010),
1312-1323.

\bibitem{gz} I. Gutman and B. Zhou, Laplacian energy of a graph, \emph{ Linear Algebra
Appl.} 414(2006), 29-37.

\bibitem{yp1} Y. Hou, Unicyclic graphs with minimal energy, {\it J.
Math. Chem.} 29(2001), 163-168.

\bibitem{yp3} Y. Hou, Bicyclic graphs with minimum energy, \emph{Linear and Multilinear
Algebra} 49(2001), 347-354.

 \bibitem{hou} Y. Hou and T. Lei, Characteristic polynomials of
skew-adjacency matrices of oriented graphs, \emph {Electron. J.
Combin.} 18(2011), $\#$p156.

\bibitem{hou1} Y. Hou, X. Shen and C. Zhang, Oriented unicyclic graphs with extremal skew energy,
Avialable at http://arXiv.org/abs/1108.6229.

\bibitem{ls} X. Li, Y. Shi and I. Gutman, \emph{Graph Energy}, Springer, New York, 2012.

\bibitem{lz1} X. Li, J. Zhang and L. Wang, On bipartite graphs with minimal
energy, \emph{Discrete Appl. Math.} 157(2009), 869-873.

\bibitem{Shen} X. Shen, Y. Hou and C. Zhang, Bicyclic digraphs with extremal skew energy,
 \emph{Electron. J. Linear Algbra} 23(2012), 340-355.

\bibitem{Tian}  G. Tian, On the skew energy of orientations of hypercubes,  \emph{ Linear Algbra
Appl.} 435(2011), 2140-2149.

\bibitem{zz} J. Zhang and B. Zhou, On bicyclic graphs with minimal energies, \emph{J.
Math. Chem.} 37(4)(2005), 423-431.

\bibitem{Zhu} J. Zhu, Oriented unicyclic graphs with the first $\lfloor \frac {n-9} 2 \rfloor$
largest skew energies, Linear Algebra Appl. 437(2012), 2630-2649.

\end{thebibliography}
\end{document}